# DENSITY ESTIMATION FOR BIASED DATA[1]

BY SAM EFROMOVICH

*The University of New Mexico*

The concept of biased data is well known and its practical applications range from social sciences and biology to economics and quality control. These observations arise when a sampling procedure chooses an observation with probability that depends on the value of the observation. This is an interesting sampling procedure because it favors some observations and neglects others. It is known that biasing does not change rates of nonparametric density estimation, but no results are available about sharp constants. This article presents asymptotic results on sharp minimax density estimation. In particular, a coefficient of difficulty is introduced that shows the relationship between sample sizes of direct and biased samples that imply the same accuracy of estimation. The notion of the restricted local minimax, where a low-frequency part of the estimated density is known, is introduced; it sheds new light on the phenomenon of nonparametric superefficiency. Results of a numerical study are presented.

**1. Introduction.** Assume that we wish to estimate the probability density $f$ of a random variable $X$. If independent direct realizations $X_1, X_2, \ldots, X_n$ of $X$ are available, then optimal solutions of the problem are well known. See the discussion in the books by Devroye and Györfi (1985), Silverman (1986) and Efromovich (1999).

In practice it may happen that drawing a direct sample from $X$ is impossible. Instead, an observation $X = x$ may be included with a relative chance proportional to a so-called biasing function $w(x)$. Then independently recorded biased observations $Y_1, Y_2, \ldots, Y_n$ have the density

(1.1) $\qquad g(y) = w(y) f(y) / \mu(f), \qquad \text{where } \mu(f) = E_f\{w(X)\}.$

The distribution of a corresponding random variable $Y$ is called a *biased* distribution and its density $g$ is given in (1.1). Given the biased sample

Received March 2001; revised July 2003.
[1]Supported by NSF Grants DMS-99-71051 and DMS-02-43606.

$Y_1, Y_2, \ldots, Y_n$ and the biasing function $w$, the problem is to estimate the underlying density $f$ with minimal mean integrated squared error over a finite interval of interest.

In what follows it is always assumed that the interval of interest is the unit interval $[0, 1]$, $0 < c_1 < w(y) < c_2 < \infty$, and $w(y)$ is Riemann integrable over the unit interval.

The following examples illustrate a few general practical settings that lead to biased data sets. (a) Let a proportion $1 - w(x)$ out of the natural frequency of $X$ be *missing*. Then the density of the observed data is (1.1). Many specific biological examples can be found in the book by Buckland, Anderson, Burnham and Laake (1993). (b) *Visibility* bias is a recognized problem in aerial survey techniques for estimating, for instance, wildlife population density. Interesting particular examples can be found in Cook and Martin (1974). (c) A sampling procedure can specifically favor larger (or smaller) observations. Two classical examples discussed in Cox (1969) are observing interevent times at some random point in time and a quality control problem of estimating fiber length distribution. (d) A rather general example is a damage model where an observation $X$ may be damaged by a destructive process depending on $X$; hence undamaged observations are biased. The interested reader can find more practical examples in the review by Patil and Rao (1977).

Let us also note that in some cases a biased sampling can be a reasonable alternative to a direct sampling. As a particular example, consider a study sponsored by the National Science Foundation (NSF) and conducted by the University of New Mexico on vegetation in the Sevilleta National Wildlife Refuge. This refuge lies 65 miles south of Albuquerque in Socorro County and includes a desert. Of particular interest is a blue gramma (*Bouteloua gracilis*), which is a native perennial that provides good grazing for wildlife and livestock; 1–2 in. tall, it grows in tufts, sod or other types of clusters of different shapes. In particular, the study is devoted to monitoring the distribution of the number of blades in a cluster, and the monitoring is based on biannual manual counts. One of the possibilities for performing the counting is to choose some areas and then count blades over these areas. This approach is manageable (after all, we are talking about desert), but experiments show that then the data are contaminated by large measurement errors. Recall that measurement errors make the problem of density estimation ill-posed and dramatically worsen accuracy of estimation [see Efromovich (1999), Chapter 3]. Thus, instead of direct counting, the area is sampled by line transects. This makes observations biased because a larger cluster has a larger probability of being intersected; however, practically negligible measurement errors make the problem dramatically simpler.

The fundamental result in the theory of biased data is from Cox (1969), where the following estimator of the cumulative distribution function was

suggested:

$$\tilde{F}(x) = \hat{\mu} n^{-1} \sum_{l=1}^{n} w^{-1}(Y_l) \mathbb{1}(Y_l \leq x), \tag{1.2}$$

where

$$\hat{\mu} = \frac{1}{n^{-1} \sum_{l=1}^{n} w^{-1}(Y_l)}. \tag{1.3}$$

For biased data sets, the Cox estimator plays the same role as the empirical distribution for direct data. An important theoretical property of the estimator is that it is a nonparametric maximum likelihood estimator. Thus, according to some general results, it is asymptotically efficient in terms of dispersion of a corresponding limit process [see the discussion in Gill, Vardi and Wellner (1988)]. Cox (1969) also suggested the first consistent kernel density estimator motivated by smoothing (1.2). Later many other density estimates, including rate optimal ones, were suggested [see the discussion in Wu and Mao (1996)]. In particular, it has been established that biasing does not affect minimax rates. Interesting results on semiparametric density estimation and their applications for moderate sample sizes were obtained by Sun and Woodroofe (1997) and Lee and Berger (2001).

On the other hand, so far no research has been conducted on sharp optimal estimation that, in particular, can shed light on Cox's long-standing question about how biasing affects density estimation. Moreover, according to Efromovich (1999, 2001), the theory of sharp estimation allows a practitioner to construct and explain the performance of data-driven estimators for small sample sizes.

This article is organized as follows. The next section presents the main theoretical results. These results and their corollaries are discussed in Section 3. Section 4 provides proofs.

**2. Minimax estimation of differentiable densities.** Let us begin by recalling a known classical result for the case of direct observations. Suppose that a random variable $X$ is distributed according to a probability density $f(x)$, $-\infty < x < \infty$, and the problem is to estimate it with a minimal mean integrated squared error over the unit interval $[0,1]$. It is assumed that $f$ is $m$-times differentiable over $[0,1]$ and belongs to a corresponding Sobolev set

$$S(m, Q) = \left\{ f(x) : f(u) = \sum_{j=0}^{\infty} \theta_j \varphi_j(u), u \in [0,1], \sum_{j=1}^{\infty} (\pi j)^{2m} \theta_j^2 \leq Q \right\}. \tag{2.1}$$

Here and in what follows, $\varphi_0(u) = 1$, $\varphi_j(u) = 2^{1/2} \cos(\pi j u)$, $j > 0$, $m$ is a positive integer number and $Q$ is a positive real number. Also define the corresponding class of densities $\mathcal{H}(m, Q) = \{f(x) : f(x) \geq 0, -\infty < x < \infty, \int_{-\infty}^{\infty} f(x)\,dx = 1, f \in S(m, Q)\}$.

For this class of densities and the case of $n$ direct observations $X_1, X_2, \ldots, X_n$, it is known that

$$(2.2) \quad \inf_{\check{f}} \sup_{f \in \mathcal{H}(m,Q)} [I_{f1} n]^{2m/(2m+1)} E_f \left\{ \int_0^1 (\check{f}(x) - f(x))^2 \, dx \right\} \geq 1 + o(1),$$

where the infimum is taken over all possible estimates $\check{f}$ based on the data set and the parameters $m$ and $Q$, and

$$(2.3) \qquad I_{f1} = Q^{-1/2m} [\pi(m+1) m^{-1} (2m+1)^{-1/2m}] \Big/ \int_0^1 f(x) \, dx.$$

Moreover, there exist data-driven estimators that attain the lower bound (2.2) [see the discussion in Efromovich (1999), Chapter 7]. Note that a typical case considered in the literature is where $[0,1]$ is the support and thus the denominator in (2.3) is equal to 1.

The approach used is called global minimax because an estimated density can be any Sobolev function. On the other hand, in practical applications an underlying density is always fixed. To bridge these two settings, Golubev (1991) suggested introducing a fixed density $f_0$, not necessarily a Sobolev one, and assuming that all possible underlying densities are uniformly close to it on the unit interval.

Namely, let $f_0$ be a density on $(-\infty, \infty)$ that is continuous and bounded below from zero on the interval $[0,1]$. No assumption about $f_0(x)$ for $x$ beyond the unit interval is made. Introduce a class of densities $D(m, Q, f_0, \rho) = \{f : \int_{-\infty}^{\infty} f(x) \, dx = 1, f(x) \geq 0, f(u) = f_0(u) + t(u), 0 \leq u \leq 1, t \in S(m, Q), \sup_{x \in [0,1]} |t(x)| < \rho\}$. Then the problem is to construct a minimax estimate for this set.

Let us present a lower bound for the local minimax approach and the case of biased data. In this case, observations $Y_1, \ldots, Y_n$ of a biased random variable $Y$ are given, the density $g(y)$ of $Y$ is defined in (1.1), and we recall that assumptions about the given biasing density $w$ are formulated below (1.1). Define

$$(2.4) \qquad\qquad\qquad\qquad I_{fw} = I_{f1}/\text{RCDB}.$$

Here RCDB is the relative coefficient of difficulty due to biasing:

$$(2.5) \quad \text{RCDB} = \int_{-\infty}^{\infty} f(x) w(x) \, dx \int_0^1 f(x) w^{-1}(x) \, dx \Big/ \int_0^1 f(x) \, dx.$$

THEOREM 1. *For any $\rho > 0$,*

$$\inf_{\check{f}} \sup_{f \in D(m,Q,f_0,\rho)} [I_{fw} n]^{2m/(2m+1)} E_f \left\{ \int_0^1 (\check{f}(x) - f(x))^2 \, dx \right\} \geq 1 + o(1),$$
(2.6)

*where the infimum is taken over all possible estimates $\check{f}$ based on the biased data set $Y_1,\ldots,Y_n$, the density $f_0$, the biasing function $w$ and parameters $m$, $Q$ and $\rho$.*

Note that (2.6) yields a corresponding global lower bound by choosing $f_0$ that is constant on the unit interval. On the other hand, under the local minimax approach, neither the set $D$ is necessarily a subset of the Sobolev set nor does $f_0$ necessarily belong to the Sobolev set.

Thus, it is absolutely natural to consider a local minimax setting where $f_0$ is more dramatically related to a class of possible underlying densities. This goal is achieved by the restricted minimax setting where $f_0$ belongs to the Sobolev set and all possible underlying densities have the same low-frequency part as $f_0$. In what follows, we refer to such $f_0$ as the anchor density, and we need the following property of $f_0$.

ASSUMPTION A (On anchor density). The anchor density $f_0(x)$ is known, positive and $m$-fold differentiable on $[0, 1]$, $f_0(u) = \sum_{j=0}^{\infty} \theta_{0j}\varphi_j(u)$, $u \in [0,1]$, and $\sum_{j=1}^{\infty}(\pi j)^{2m}\theta_{0j}^2 = Q$. Also, there exists a sequence $k_s \to \infty$, $s \to \infty$, such that $\sum_{j>k_s} j^{2m}\theta_{0j}^2 > C_1 k_s^{-C_2}$ for some positive $C_1$ and $C_2$.

Let us comment on the two parts of the assumption made. The first part implies that the anchor density is a particular density from the Sobolev set (2.1) studied under the global minimax setting. Note that either the statistician may know $Q$ and then choose a corresponding anchor density, or the statistician may choose an anchor density and then calculate $Q$. The second part (the part about the existence of $k_s$) assumes that the anchor density is not too smooth and thus it is a typical density from the Sobolev set. For instance, let us check that the assumption about $k_s$ holds whenever $\sum_{j>0} j^{2m+\alpha}\theta_{0j}^2 = \infty$ for some $\alpha > 0$. Indeed, if no sequence $k_s$ exists for a particular $C_2 = 1 + 2\alpha$, then $\sum_{j>k} j^{2m}\theta_{0j}^2 < Ck^{-1-2\alpha}$, $C < \infty$, for all sufficiently large $k$. The last inequality implies $\theta_{0j}^2 < Cj^{-2m-1-2\alpha}$ and thus we get the inequality $\sum_{j>0} j^{2m+\alpha}\theta_{0j}^2 < \infty$ that contradicts the assumption made.

The second part of the assumption yields that there always exists an increasing to infinity integer-valued sequence $J_n$ such that $-1 \leq J_n < \ln(n)$ and $\sum_{j>J_n} j^{2m}\theta_{0j}^2 > C_1(J_n + 2)^{-C_2}$. From now on this sequence is assumed to be fixed. Then we introduce the sequence of low-frequency parts of $f_0$,

$$(2.7) \qquad f_{0J_n}(u) = \sum_{j=0}^{J_n} \theta_{0j}\varphi_j(u),$$

and define the vanishing sequence

$$
(2.8) \qquad q_n = 1 - \frac{\sum_{j=1}^{J_n} j^{2m} \theta_{0j}^2}{\sum_{j=1}^{\infty} j^{2m} \theta_{0j}^2}.
$$

In what follows it is assumed that $\sum_{j=a}^{b} c_j = 0$ whenever $b < a$.

Now we are in position to define the restricted density set:

$$
(2.9) \qquad \mathcal{H}(m, Q, f_0, J_n) = \bigg\{ f : f(u) = f_{0 J_n}(u) + \sum_{j > J_n} \theta_j \varphi_j(u), u \in [0,1], \\
f \in S(m, Q), f(x) \geq 0, \int_{-\infty}^{\infty} f(x)\, dx = 1 \bigg\}.
$$

It is easy to see that $\max_{x \in [0,1]} |f(x) - f_0(x)| = o(1)$ uniformly over $f \in \mathcal{H}(m, Q, f_0, J_n)$; thus the restricted approach is also local around the anchor density $f_0$. Also, we may set $J_n = -1$ and then the restricted approach becomes global.

THEOREM 2. *Let Assumption* A *hold. Then*

$$
(2.10) \qquad \inf_{\check{f}} \sup_{f \in \mathcal{H}(m,Q,f_0,J_n)} [I_{fw} q_n^{-1/2m} n]^{2m/(2m+1)} E_f \bigg\{ \int_0^1 (\check{f}(x) - f(x))^2\, dx \bigg\} \\
\geq 1 + o(1),
$$

*where the infimum is taken over all possible estimates $\check{f}$ based on the biased data set $Y_1, \ldots, Y_n$, the anchor density $f_0$, the biasing function $w$, the parameters $m$, $Q$ and the sequence $J_n$.*

The lower bounds (2.6) and (2.10) are attained by the Efromovich–Pinsker adaptive estimator, which is a blockwise shrinkage estimator defined as follows. We divide the set of natural numbers into a sequence of nonoverlapping blocks $G_k$, $k = 1, 2, \ldots$. Then the estimate is

$$
(2.11) \qquad \hat{f}(x) = \sum_{k=1}^{K} \bigg[ 1 - \frac{\hat{d} n^{-1}}{|G_k|^{-1} \sum_{j \in G_k} \hat{\theta}_j^2} \bigg] \\
\times \mathbb{1}\bigg( |G_k|^{-1} \sum_{j \in G_k} \hat{\theta}_j^2 > (1 + t_k) \hat{d} n^{-1} \bigg) \sum_{j \in G_k} \hat{\theta}_j \varphi_j(x),
$$

where $|G_k|$ denotes the cardinality of $G_k$, $\mathbb{1}(\cdot)$ is the indicator,

$$
(2.12) \qquad \hat{\theta}_j = \hat{\mu} n^{-1} \sum_{l=1}^{n} \mathbb{1}(0 \leq Y_l \leq 1) \varphi_j(Y_l) w^{-1}(Y_l)
$$

is the Cox sample mean estimate of Fourier coefficients and

$$\hat{d} = \hat{\mu}^2 n^{-1} \sum_{l=1}^{n} \mathbb{1}(0 \leq Y_l \leq 1) w^{-2}(Y_l). \tag{2.13}$$

A wide variety of blocks $\{G_k\}$ and thresholds $\{t_k\}$ implies sharp minimaxity [see the discussion in Efromovich (1985, 1999, 2000)]. As an example, we set $|G_k| = k^2$, $t_k = 1/\ln(k+1)$ and $K = \lfloor n^{1/9} \ln(n) \rfloor$.

THEOREM 3. *The Efromovich–Pinsker estimator satisfies*

$$\sup_{f \in \mathcal{H}(m,Q)} [I_{fw} n]^{2m/(2m+1)} E_f \left\{ \int_0^1 (\hat{f}(u) - f(u))^2 \, du \right\} \tag{2.14}$$
$$= 1 + o(1),$$

*and if Assumption* A *holds, then*

$$\sup_{f \in \mathcal{H}(m,Q,f_0,J_n)} [I_{fw} q_n^{-1/2m} n]^{2m/(2m+1)} E_f \left\{ \int_0^1 (\hat{f}(u) - f(u))^2 \, du \right\} \tag{2.15}$$
$$= 1 + o(1).$$

## 3. Discussion.

3.1. *The minimax approaches.* It may be convenient to think about the minimax approaches in terms of concepts of game theory. We may think that nature (Player I) chooses a density that makes its estimation most difficult for the statistician (Player II). Then the main difference between the three minimax approaches introduced in Section 2 is in the information available to the statistician. Under the global approach, the statistician knows that nature chooses a density from a given Sobolev set. Under the local approach, the statistician knows that nature chooses a density which is uniformly close to a given density.

Under the restricted approach the statistician knows dramatically more about nature's choice. The statistician has the same information as in the global game and additionally knows the low-frequency part of nature's choice. The latter also makes the game local because an underlying density is uniformly close to its low-frequency part. In other words, the restricted game is about estimating the high-frequency part of an underlying Sobolev density.

A minimax data-driven estimator (i.e., the estimator based only on data and the biasing function) should perform not worse than the statistician playing the minimax game. Thus the restricted minimax game is more challenging for a data-driven estimation. On the other hand, the restricted game is more rewarding because the rate of mean integrated squared error convergence is always faster than the global or local minimax rate $n^{-2m/(2m+1)}$.

It has long been a tradition in the nonparametric literature to study global and/or local minimax approaches. This explains the familiar slogan "...if we are prepared to assume that the unknown density has $m$ derivatives, then ...the optimal mean integrated squared error is of order $n^{-2m/(2m+1)}$...." The citation is from Silverman [(1986), page 70]. The results of Section 2 show that this classical rate is optimal only if data-driven estimates are compared with the statistician playing global or local minimax games, that is, with the less informed statistician. Faster rates can be obtained by matching the restricted minimax game.

Let us make one more remark about the minimax approaches. It is possible to change the setting a bit and to assume that the underlying density is known beyond the interval of interest. This makes the average value $\int_0^1 f(x)\,dx$ known, but this fact does not affect the asymptotics.

3.2. *Practical implications of the minimax approaches.* At first glance, because restricted minimax implies faster asymptotic rates, there is no reason to study global and local minimax approaches.

Interestingly, small data sets justify the study of these classical minimax approaches. It was shown by Efromovich (1999) that, for small sample sizes (up to several hundred), the problem of nonparametric density estimation is equivalent to the problem of estimating a low-frequency part of the underlying density. The reader can check this assertion using the software in Efromovich [(1999), Chapter 3]. As a result, the restricted minimax approach with nonnegative $J_n$ is not applicable for small sample sizes, because it assumes that a low-frequency part of the underlying density is given. (This is also the reason behind the construction of $J_n$ that allows us to apply no restrictions on the underlying density for small $n$.)

The situation changes for moderate and large sample sizes like the ones studied in the wavelet literature. For these sample sizes, knowing or not knowing a low-frequency part of the density has no significant effect on the estimation, and thus the restricted minimax is absolutely appropriate. Again, the interested reader can use the software to check this assertion.

We may conclude that each minimax approach has its own practical applications.

3.3. *Restricted minimax and nonparametric superefficiency.* The phenomenon of parametric superefficiency is well known. A famous example is the Hodges superefficient estimator that, for normal observations, allows us to improve a sample mean estimator (efficient estimator) at any given point. This is an interesting theoretical phenomenon; on the other hand, superefficient estimators are typically not used, because the set of super-efficiency has Lebesgue measure zero and estimation at other points may worsen. See the discussion in Ibragimov and Khasminskii [(1981), Section 2.13].

By contrast, it was shown in Brown, Low and Zhao (1997) that, in nonparametric problems, every curve can be a point of superefficiency. Their main result, "translated" into our density estimation setting, is that for any $f \in \mathcal{H}(m,Q)$ there exists an estimator $\check{f}_n$ such that

$$(3.1) \qquad n^{2m/(2m+1)} E_f\left\{\int_0^1 (\check{f}_n(x) - f(x))^2\, dx\right\} = o(1).$$

This result implies a better rate than the classical $n^{-2m/(2m+1)}$. This explains why Brown, Low and Zhao (1997) refer to (3.1) as the nonparametric superefficiency.

On the other hand, (3.1) is in agreement with the restricted minimax rates. Let us also note that if an underlying density is parametric (it has a finite number of nonzero Fourier coefficients), then the Efromovich–Pinsker estimator implies the parametric rate $n^{-1}$ of the mean integrated squared error convergence [see Efromovich (1985)]. This indicates the range of possible nonparametric rates.

3.4. *Restricted minimax and oracles.* An oracle is an estimator that is based on both data and the underlying density. The oracle approach means a data-driven estimation that mimics the oracle performance [see the discussion in Efromovich (1999)]. The restricted minimax bridges classical minimax approaches (where underlying densities belong to function spaces) and oracle approaches (where the underlying density is given) by assuming that the underlying density belongs to a function space and its low-frequency part is given.

3.5. *Average risk.* One of the long-standing problems in the sharp estimation literature is to find an estimate of the density $g(y)$ of direct observations $Y_1, Y_2, \ldots, Y_n$ that minimizes an average risk with the averaging function $a(y)$. In other words, this estimate should minimize the average risk

$$(3.2) \qquad \mathrm{AR} = E_g\left\{\int_0^1 a(y)(\hat{g}(y) - g(y))^2\, dy\right\}.$$

Let us assume that the averaging function satisfies $0 < C_* \le a(y) \le C^* < \infty$. Then it is easy to see that if we define

$$(3.3) \qquad \begin{aligned} w(y) &= 1/\sqrt{a(y)}, \\ f(y) &= \mu(f) w^{-1}(y) g(y), \\ \hat{f}(y) &= \mu(f) w^{-1}(y) \hat{g}(y), \end{aligned}$$

then

$$(3.4) \qquad \mathrm{AR} = \mu^{-2}(f) E_f \left\{ \int_0^1 (\hat{f}(x) - f(x))^2 \, dx \right\}.$$

Here $f$ and $g$ can be thought of as the underlying and the biased densities.

Particular examples of using this equivalence for finding sharp asymptotics are presented in Efromovich (2004).

3.6. *Naive estimation.* Using a naive estimator $\tilde{f}_n(x) = \tilde{g}_n(x)\hat{\mu}w^{-1}(x)$ is a popular and intuitively clear idea [see the discussion in Wu and Mao (1996) and Wu (1997)]. Here $\tilde{g}_n$ is an estimator of the density $g$ of biased observations $Y_1, \ldots, Y_n$. Section 2 implies that smoothness of the biasing function plays a crucial role in the accuracy of the naive estimator. The naive estimator is rate inadmissible whenever the biasing function is not as smooth as the underlying density $f$. On the other hand, specific examples where naive estimation is sharp minimax can be found in Efromovich (2004). We may conclude that because smoothness of $f$ is typically unknown, it is better to avoid the use of naive estimation.

It is also important to note that smoothness of the biasing function does not affect optimal estimation. Thus, even if the biased distribution is discontinuous or not differentiable, the quality of sharp minimax estimation of the underlying density $f$ is defined only by its own smoothness. The only functional of $w$ that affects the estimation is the coefficient of difficulty due to biasing, discussed in the next section.

3.7. *Coefficient of difficulty.* An interesting theoretical outcome of Section 2 is that, for a biased sample of size $n$, the same precision of estimation is achievable by a direct sample of size $n' = n/\mathrm{RCDB}$, where RCDB is the relative coefficient of difficulty due to biasing defined in (2.5). Recall that the notion of the coefficient of difficulty was introduced in Efromovich (1999).

Thought-provoking examples in Cox (1969) indicate that biasing can always improve or worsen the estimation. Translated into the nonparametric setting considered, this would imply that biasing can always increase or decrease the RCDB.

Let us present an example where RCDB is always greater than 1, that is, the example where biasing always worsens the density estimation. According to the Cauchy–Schwarz inequality,

$$(3.5) \quad \left[ \int_{-\infty}^{\infty} \mathbb{1}(0 \le x \le 1) f(x) \, dx \right]^2 \le \int_{-\infty}^{\infty} f(x) w(x) \, dx \int_0^1 f(x) w^{-1}(x) \, dx$$

with equality iff $w(x) = c\mathbb{1}(0 \le x \le 1)$, $c > 0$, almost sure with respect to $f$. Thus, if $[0,1]$ is the support of $X$, then any biasing yields $\mathrm{RCDB} > 1$, that is, biasing always worsens the density estimation. Otherwise, $\mathrm{RCDB} \ge$ $\int_0^1 f(x)\,dx$ and, similarly to examples in Cox (1969), biasing can improve or worsen the estimation.

This is a useful conclusion for practitioners because, as we shall see in Section 3.9, the asymptotic RCDB can be used for the analysis of small data sets.

3.8. *Versatility of the Efromovich–Pinsker estimator.* It is well known that many functionals of this estimator (including derivatives and integrals) are optimal estimators of the corresponding functionals, that is, the estimator is versatile [see the discussion in Efromovich (1999), Chapter 7]. It is possible to establish that the same conclusion holds for the case of biased data; and results will be published elsewhere.

3.9. *Adaptive estimation*: *from asymptotic to small sample sizes.* Asymptotic results presented in Section 2 show that the Efromovich–Pinsker series estimator is asymptotically minimax. This justifies the use of software developed in Efromovich [(1999), Chapter 3] for small data sets. The software, which contains both a generator of biased data sets and the adaptive estimator for small biased data sets, is available over the Worldwide Web [see the instructions on how to download and use it in Efromovich (1999), Appendix B or e-mail the author].

Using this software, let us shed light on the nature of a biased sampling and then comment on the possibility of using the coefficient of difficulty RCDB for small data sets.

Figure 1 presents a particular biased data set of size $n = 25$ shown by letters Y. The underlying density $f$ is the Normal density shown by the solid line and defined in Efromovich [(1999), page 18]. The sample is biased by the biasing function $w(y) = 0.1 + 0.9y$ shown by the long-dashed line, that is, the data may be referred to as length biased. The right-skewed data set clearly exhibits the effect of this biasing. To exhibit the structure of the data set, the short-dashed line shows us its estimated density [the adaptive estimate of Efromovich (1999) is used]. This is what the statistician might see if the biased nature of the data were ignored or unknown. Note how the skewed density of Y's differs from the symmetric Normal density.

The dotted line shows the suggested adaptive biased-data density estimate. By taking into account the biasing, the estimate correctly restored the symmetric about 0.5 shape of the underlying density. It also removed the heavy left tail of the estimated density of Y's created by the three smallest length-biased observations.

This particular simulation together with the discussion in the Section 3.7 raises the following question. Suppose that $n'$ is a sample size that implies a reasonable estimation of an underlying density based on direct observations. Then what is the corresponding sample size $n$ for a biased data set that

## Biased Data

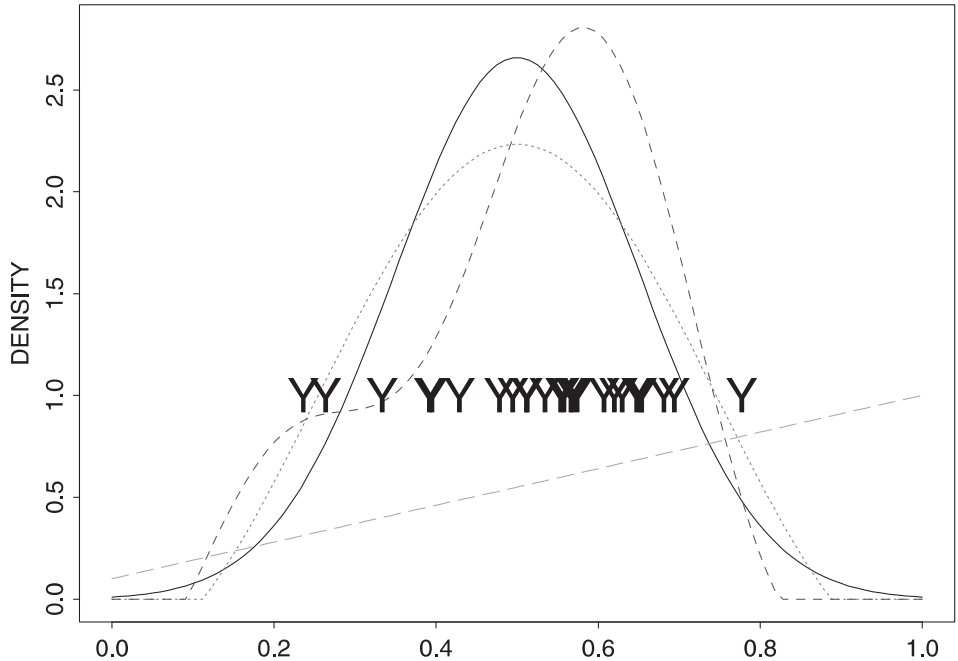

Fig. 1. *Analysis of biased data. The solid and long-dashed lines show the underlying Normal density $f$ and the biasing function $w(y) = 0.1 + 0.9y$, respectively. A simulated biased data set of size $n = 25$ is shown by letters Y. The dotted line shows the estimate of $f$ (the estimate based on Y's and the biasing function). The short-dashed line shows the estimate of $g$, that is, of the density of the biased observations Y.*

implies a similar precision of estimation in terms of mean integrated squared error? According to the asymptotic results of Section 2, $n'$ times RCDB should be the answer, but can this asymptotic rule be used for small sample sizes?

Let us begin the discussion with a particular simulation and then complement it with an intensive Monte Carlo study.

Figure 2 explains the problem explored. The underlying density is the monotone one shown by the solid line and defined in Efromovich [(1999), page 18]. A direct sample of size $n' = 25$ from this density is shown by X's. Note that the sample correctly represents the underlying density, and this also can be seen from the estimate of the density shown by the short-dashed line.

A biased sample of size $n = 44$ from the same density is shown by Y's. The utilized biasing function is $w(y) = 1 - 0.95y$, and this implies RCDB $= 1.74$ and the above-mentioned sample size $n = n' \times$ RCDB $= 44$. If the asymptotic theory holds for these small sample sizes, then the density estimation based on 25 direct and 44 biased data values should be similar in terms of mean integrated squared error. For the particular samples the estimates for direct and biased samples are shown by the short-dashed and dotted lines, respectively. The short-dashed line better exhibits the underlying density, and it may look like the 1.74-fold increase in the sample size is not large enough to compensate for the biasing. On the other hand, let us recall that the samples are independent and another simulation may change the outcome.

Let us repeat this particular simulation 500 times, calculate corresponding integrated squared errors (ISEs) and then analyze them. The results are presented in Figure 3. The top diagram shows by character 1 ISEs for direct data sets and by 2 for biased data sets. Let us repeat that all samples are independent. Clearly a majority of ISEs are relatively small but there is a thin right tail in their distribution. Thus, we show densities of ISEs over two subintervals: ISEs that are at most 0.11 and larger ISEs. The value 0.11 is the mean (up to the rounded second digit) of both sets of ISEs (here we have an ideal outcome in terms of the empirical mean integrated squared errors).

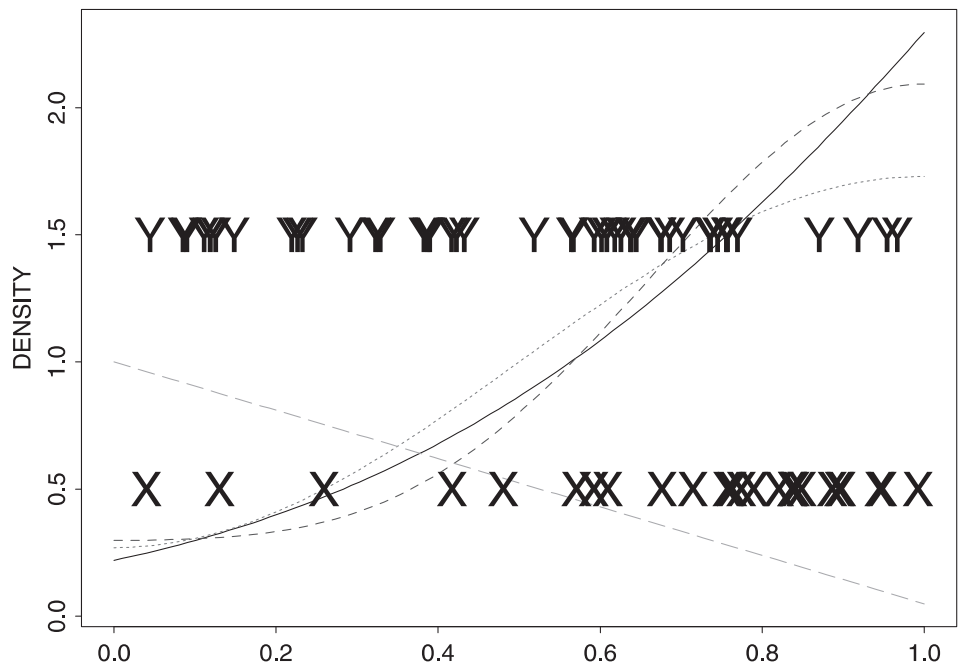

FIG. 2. *Analysis of biased and direct data. The solid and long-dashed lines show the underlying monotone density $f$ and the biasing function $w(y) = 1 - 0.95y$, respectively. The corresponding $\mathrm{RCDB} = 1.74$. A simulated direct data set of size $n' = 25$ from the underlying monotone density is shown by letters* X. *A biased data set of adjusted size $n = n' \times \mathrm{RCDB} = 44$ is shown by letters* Y. *The dotted line shows the biased-data density estimate (the estimate is based on* Y*'s and the biasing function) and the short-dashed line shows the density estimate of* X.

The densities for these two groups of ISEs are shown in the middle and bottom diagrams. As we see, the distributions are practically identical, and while there is no asymptotic theory to support this outcome, it is an interesting empirical observation.

Figure 4 shows an outcome of a similar study only with the underlying density and the biasing function utilized in Figure 1. The main parameters are presented in the caption, and here let us stress only the relatively small RCDB = 1.07. As we see, the outcome is very similar. The sample means

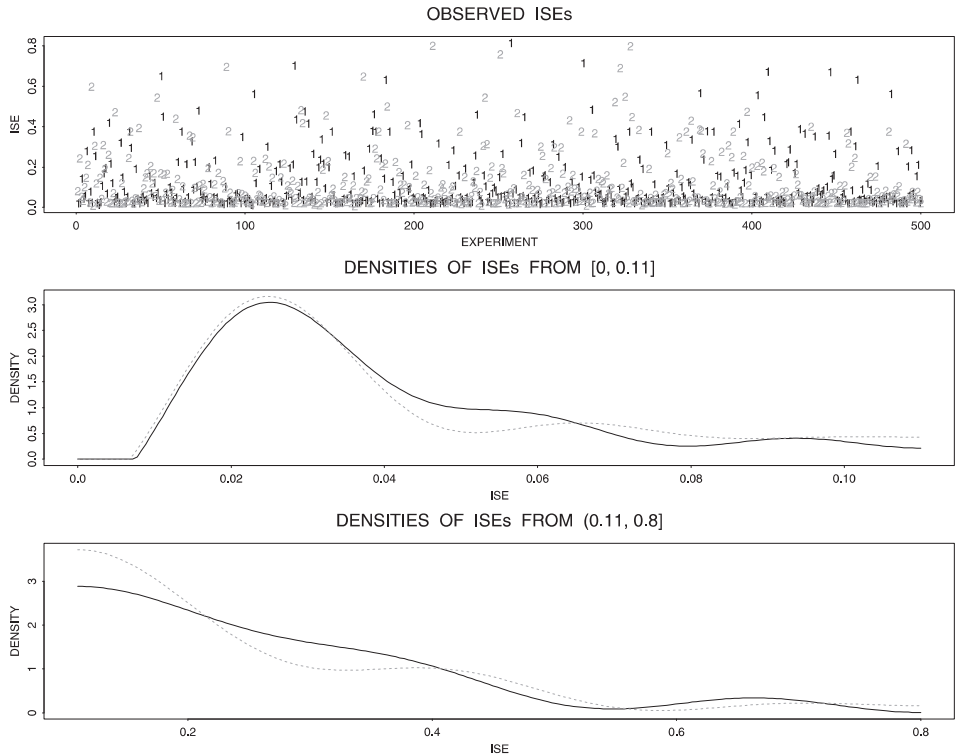

Fig. 3. *Results of* 500 *Monte Carlo simulations identical to the one shown in Figure* 2. *Characters* 1 *and* 2 *in the top diagram show ISEs of the estimates based on* 25 *direct and* 44 *biased observations, respectively. The sample means are identical (up to the second digit) and equal to* 0.11. *The solid and dotted lines in the bottom diagrams show the densities of ISEs for the direct and biased samples, respectively. The two bottom diagrams show the densities for ISEs that are at most the sample mean* 0.11 *and larger than the sample mean, respectively. From the totals of* 500, *there are* 348 *and* 357 *ISEs that are at most* 0.11 *for the direct and biased samples, respectively.*

DENSITY ESTIMATION FOR BIASED DATA 15

are a bit different (they are 0.09 and 0.10 for the direct and biased samples, resp.), but it is clear that the difference is primarily due to the tails. Repeated simulations show that this is indeed the case.

The numerical study supports the possibility of using RCDB as a measure of difficulty due to biasing. The interested reader can find a different numerical study, which implies a similar outcome for a wider variety of densities and biasing functions, in Efromovich (2004).

**4. Proofs.** Recall that $\varphi_0(x) = 1$, $\varphi_j(x) = \sqrt{2}\cos(\pi j x)$, $j \geq 1$, and $\lfloor x \rfloor$ is the rounded down $x$.

PROOF OF THEOREM 1. This proof will be also used to verify Theorem 2, and this explains some steps and comments not directly related to the proof.

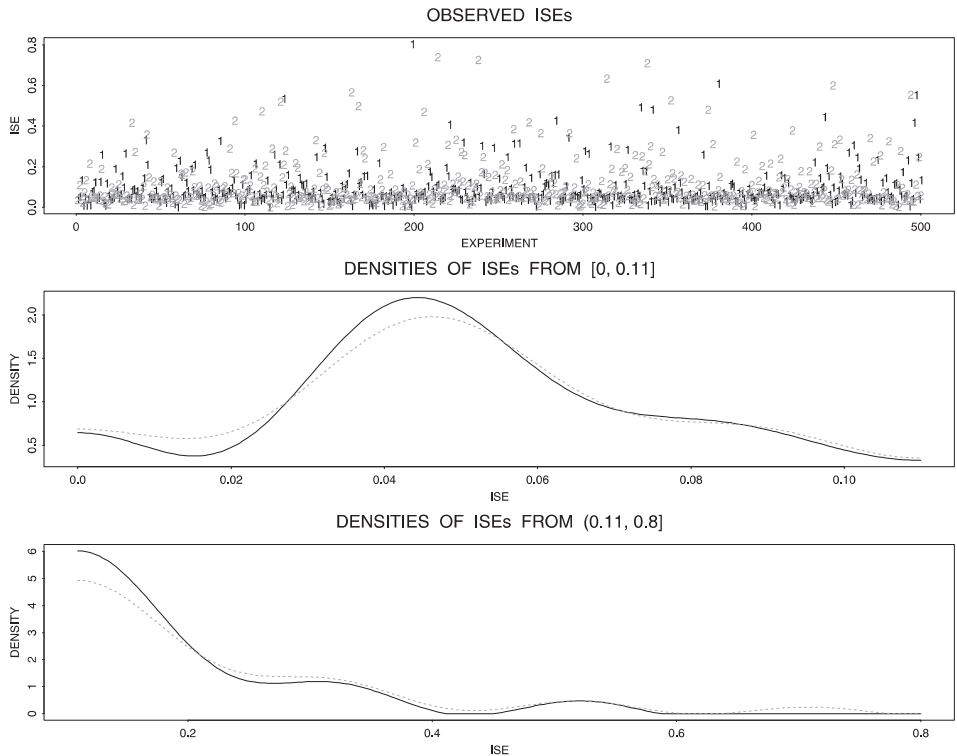

FIG. 4. *A numerical study similar to the one shown in Figure* 3, *only here the density and the biasing function of Figure* 1 *are used.* RCDB $= 1.07$, *and this implies* 25 *direct and* 27 *biased observations. Mean ISEs are* 0.09 *and* 0.10, *respectively. From the totals of* 500, *there are* 383 *and* 362 *ISEs that are at most* 0.11 *for the direct and biased samples, respectively.*

We begin by dividing the unit interval into $s$ subintervals where the biasing function $w$ and the density $f_0$ are approximated by simple functions. This allows us to obtain relatively simple lower minimax bounds for each subinterval.

Set $s = 1 + \lfloor \ln(\ln(n+20)) \rfloor$ and define

$$\mathcal{H}_s = \left\{ f : f(x) = f_0(x) + \left[ \sum_{k=0}^{s-1} f_k(x) - \sum_{k=0}^{s-1} \int_0^1 f_k(u)\, du \right] \mathbb{1}(0 \leq x \leq 1), \right.$$

$$\left. f_k(x) \in \mathcal{H}_{sk}, f \geq 0 \right\}.$$

Here the function classes $\mathcal{H}_{sk}$ are defined as follows. Let $\phi(x) = \phi(n, x)$ be a sequence of flattop nonnegative kernels defined on a real line such that, for a given $n$, the kernel is zero beyond $(0,1)$, it is $m$-fold continuously differentiable on $(-\infty, \infty)$, $0 \leq \phi(x) \leq 1$, $\phi(x) = 1$ for $2(\ln(n))^{-2} \leq x \leq 1 - 2(\ln(n))^{-2}$ and $|\phi^{(m)}| \leq C(\ln(n))^{2m}$. For instance, such a kernel may be constructed using so-called mollifiers discussed in Efromovich [(1999), Chapter 7]. Then set $\phi_{sk}(x) = \phi(sx - k)$. For the $k$th subinterval, $0 \leq k \leq s-1$, define $\varphi_{skj}(x) = \sqrt{s}\varphi_j(sx - k)$, $f_{[k]}(x) = \sum_{j=\lfloor J(k)/\ln(n) \rfloor}^{J(k)} \nu_{skj} \varphi_{skj}(x)$, $f_{(k)}(x) = f_{[k]}(x)\phi_{sk}(x)$, $J(k) = 2\lfloor [n(2m+1)(m+1)s^{-2m} Q_{sk}(2m(2\pi)^{2m})^{-1}]^{1/(2m+1)} \rfloor$, $Q_{sk} = Q(1-1/s)(\overline{I_s^{-1}} I_{sk})^{-1}$, $I_{sk} = \mu^{-1}(f_0) w(k/s)/f_0(k/s)$, and $\overline{I_s^{-1}} = \sum_{k=0}^{s-1}(1/I_{sk})$. Then we define the subclasses

$$\mathcal{H}_{sk} = \left\{ f : f(x) = f_{(k)}(x), \sum_{j=\lfloor J(k)/\ln(n) \rfloor}^{J(k)} (\pi j)^{2m} \nu_{skj}^2 \leq s^{-2m} Q_{sk}, \right.$$

$$\left. |f_{[k]}(x)|^2 \leq s^3 \ln(n) J(k) n^{-1} \right\}.$$

Let us verify that, for sufficiently large $n$, this set of densities is a subset of the studied class $D$.

First, the definition of the flattop kernel implies that $f - f_0$ is $m$-fold continuously differentiable over $[0,1]$. Second, let us verify that for $f \in \mathcal{H}_s$ the difference $f - f_0$ belongs to $S(m, Q)$. By the Leibniz rule, $\int_0^1 [(f_{[k]}(x)\phi_{sk}(x))^{(m)}]^2\, dx = \int_0^1 [\sum_{l=0}^m \mathbf{C}_l^m f_{[k]}^{(m-l)}(x) \phi_{sk}^{(l)}(x)]^2\, dx$, where $\mathbf{C}_l^m = m!/((m-l)!l!)$. Recall that $\max_{0 \leq l \leq m} \int_0^1 (\phi_{sk}^{(l)}(x))^2\, dx < C(s(\ln(n))^2)^{2m}$ and, for $0 < l \leq m$,

$$|f_{[k]}^{(m-l)}(x)|^2 = \left| \sum_{j=\lfloor J(k)/\ln(n) \rfloor}^{J(k)} \nu_{skj} \varphi_{skj}^{(m-l)}(x) \right|^2$$

$$(4.1) \qquad \leq C s^{2(m-l)+1} \left( \sum_{j=\lfloor J(k)/\ln(n) \rfloor}^{J(k)} j^{2m} \nu_{skj}^2 \right) \left( \sum_{j=\lfloor J(k)/\ln(n) \rfloor}^{J(k)} j^{-2l} \right)$$

$$= o(1)(J(k))^{-1/2},$$

where the Cauchy–Schwarz inequality was used in the middle line. Also,

$$(4.2) \quad \int_0^1 [f_{[k]}^{(m)}(x)\phi_{sk}(x)]^2\, dx \leq \int_{k/s}^{(k+1)/s} (f_{[k]}^{(m)}(x))^2\, dx \leq Q_{sk},$$

and recall that $\sum_{k=0}^{s-1} Q_{sk} = Q(1-s^{-1})$. These results imply $f - f_0 \in S(m, Q(1-s^{-2}))$ for $f \in \mathcal{H}_s$ and large $n$.

Now denote

$$\hat{f} = f_0 + \tilde{f} \quad \text{and} \quad \delta_s = \sum_{k=0}^{s-1} \int_0^1 f_{(k)}(u)\, du,$$

and note that for $f \in \mathcal{H}_s$ and any $\gamma > 0$,

$$\int_{k/s}^{(k+1)/s} (\hat{f}(x) - f(x))^2\, dx$$

$$= \int_{k/s}^{(k+1)/s} (\tilde{f}(x) - f_{(k)}(x) + \delta_s)^2\, dx$$

$$\geq (1-\gamma) \int_{k/s}^{(k+1)/s} (\tilde{f}(x) - f_{[k]}(x))^2\, dx$$

$$- \gamma^{-1} \int_{k/s}^{(k+1)/s} [f_{[k]}(x)(1 - \phi_{sk}(x)) + \delta_s]^2\, dx$$

$$\geq (1-\gamma) \int_{k/s}^{(k+1)/s} (\tilde{f}(x) - f_{[k]}(x))^2\, dx + o(1)\gamma^{-1}(\ln(n))^{-1/2} n^{-2m/(2m+1)}.$$

Set $\gamma = s^{-1}$ and, using the above-obtained relationship, we get

$$\sup_{f \in D(m,Q,f_0,\rho)} E\left\{\int_0^1 (\hat{f}(x) - f(x))^2\, dx\right\}$$

$$\geq \sup_{f \in \mathcal{H}_s} E\left\{\int_0^1 (\hat{f}(x) - f(x))^2\, dx\right\}$$

$$= \sup_{f \in \mathcal{H}_s} \sum_{k=0}^{s-1} E\left\{\int_{k/s}^{(k+1)/s} (\hat{f}(x) - f(x))^2\, dx\right\}$$

$$\geq (1-s^{-1}) \sum_{k=0}^{s-1} \sup_{f \in H_{sk}} \sum_{j=\lfloor J(k)/\ln(n)\rfloor}^{J(k)} E\{(\tilde{\nu}_{skj} - \nu_{skj})^2\} + o(1) n^{-2m/(2m+1)}$$

$$= (1-s^{-1}) \sum_{k=0}^{s-1} R_k + o(1) n^{-2m/(2m+1)},$$

where $\tilde{\nu}_{skj} = \int_{k/s}^{(k+1)/s} \tilde{f}(x)\varphi_{skj}(x)\,dx$.

To estimate $R_k$, following the proof of Theorem 1 in Efromovich (1989), we make two additional steps. The first one is to verify that if $\zeta_{skj}$ are independent normal random variables with zero mean and variance $(1-\gamma)\nu_{skj}^2$, where here $\gamma = \gamma_n$ tends to zero as slowly as desired, then a stochastic process $f^*(x)$, defined as the studied $f(x)$ but with $\zeta_{skj}$ in place of $\nu_{skj}$, satisfies the relationship

$$(4.3) \qquad P(f^*(x) \in \mathcal{H}(m,Q)) = 1 + o(1),$$

and if additionally $\nu_{skj}^2 \leq sn^{-1}$, then a similarly defined stochastic process $f_{[k]}^*$ satisfies

$$(4.4) \qquad P\left(\sup_{x \in [0,1]} |f_{[k]}^*(x)|^2 \leq s^3 \ln(n) J(k) n^{-1}\right) = 1 + o(1).$$

The second step is to compute for $f \in \mathcal{H}_s$ the classical parametric Fisher information

$$(4.5) \qquad I_{skj} = E_{f_0}\{[\partial \ln(f(Y)w(Y)/\mu(f))/\partial \nu_{skj}]^2\}.$$

Relationship (4.3) follows from (A.18) in Pinsker (1980). Also, for $\nu_{skj}^2 \leq sn^{-1}$, the inequality

$$\sum_{j=\lfloor J(k)/\ln(n) \rfloor}^{J(k)} \sup_x [\nu_{skj}\varphi_{skj}(x)]^2 \ln(sJ(k)) \leq Cs^2 n^{-1} J(k) \ln(n)$$

holds and this together with Theorem 6.2.2 in Kahane (1985) yields (4.4).

Now we are in a position to calculate the Fisher information (4.5). To simplify notation, let us denote $\nu_{skj} = \theta$ and the corresponding density by $f_\theta$. Write

$$\frac{\partial \ln(g_\theta(u))}{\partial \theta} = \frac{\partial \ln(w(u)f_\theta(u)/\mu(f_\theta))}{\partial \theta}$$
$$= \frac{f'_\theta(u)\mu(f_\theta) - \mu'(f_\theta)f_\theta(u)}{f_\theta(u)\mu(f_\theta)}.$$

Here $f'_\theta(u) = \partial f_\theta(u)/\partial \theta$ and $\mu'(f_\theta) = \partial \mu(f_\theta)/\partial \theta$. This implies that

$$\left[\frac{\partial \ln(g_\theta(u))}{\partial \theta}\right]^2 g_\theta(u)$$
$$= \frac{w(u)[(f'_\theta(u))^2\mu^2(f_\theta) - 2\mu(f_\theta)\mu'(f_\theta)f_\theta(u)f'_\theta(u) + (\mu'(f_\theta)f_\theta(u))^2]}{\mu^3(f_\theta)f_\theta(u)}$$
$$(4.6) \qquad = w(u)\mu^{-3}(f_\theta)[\mu^2(f_\theta)(f'_\theta(u))^2 f_\theta^{-1}(u)$$

$$-2\mu(f_\theta)\mu'(f_\theta)f'_\theta(u) + (\mu'(f_\theta))^2 f_\theta(u)]$$
$$= w(u)\mu^{-1}(f_\theta)(f'_\theta(u))^2 f_\theta^{-1}(u)$$
$$-2\mu^{-2}(f_\theta)w(u)\mu'(f_\theta)f'_\theta(u) + \mu^{-3}(f_\theta)w(u)(\mu'(f_\theta))^2 f_\theta(u)$$
$$= T_1(u) + T_2(u) + T_3(u).$$

Recall that $\phi_{ks}(x)$ is supported on $[k/s, (k+1)/s]$ and we are estimating three components of the Fisher information that correspond to the three terms on the right-hand side of (4.6). Write

$$T_1 = \int_{k/s}^{(k+1)/s} w(u)\mu^{-1}(f_0) f_0^{-1}(u)(1+o(1))$$
$$\times \left[\varphi_{skj}(u)\phi_{sk}(u) - \int_{k/s}^{(k+1)/s} \varphi_{skj}(z)\phi_{sk}(z)\,dz\right]^2 du.$$

To estimate $T_1$ we use the the following three relationships. Write

$$\int_{k/s}^{(k+1)/s} [\varphi_{skj}(x)\phi_{sk}(x)]^2 \, dx = 1 + \int_{k/s}^{(k+1)/s} \varphi_{skj}^2(x)(\phi_{sk}^2(x) - 1)\,dx$$

and then, recalling that $\phi_{sk}(x)$ is the special flattop kernel,

$$\left|\int_{k/s}^{(k+1)/s} \varphi_{skj}^2(x)(\phi_{sk}^2(x) - 1)\,dx\right| = o(1)(\ln(n))^{-1}.$$

Similarly,

(4.7)
$$\left|\int_{k/s}^{(k+1)/s} \varphi_{skj}(x)\phi_{sk}(x)\,dx\right| = \left|\int_{k/s}^{(k+1)/s} \varphi_{skj}(x)[\phi_{sk}(x) - 1]\,dx\right|$$
$$= o(1)(\ln(n))^{-1}.$$

Also, using the assumptions about the biasing function and the anchor density, we obtain that

$$T_1 = \mu^{-1}(f_0)w(k/s)f_0^{-1}(k/s)(1+o(1)) = I_{sk}(1+o(1)).$$

Using (4.7), the second component $T_2$ of the Fisher information can be estimated as

$$T_2 = -2\mu^{-2}(f_0)\mu'(f_0)$$
$$\times \int_0^1 f_0(u)\left[\varphi_{skj}(u)\phi_{sk}(u) - \int_{k/s}^{(k+1)/s} \varphi_{skj}(z)\phi_{sk}(z)\,dz\right]du(1+o(1))$$
$$= o(1).$$

To estimate $T_3$ we write

$$\mu'(f_\theta) = \int_0^1 f'_\theta(u) w(u)\, du$$

$$= \int_0^1 \left[ \varphi_{skj}(u)\phi_{sk}(u) - \int_{k/s}^{(k+1)/s} \varphi_{skj}(z)\phi_{sk}(z)\, dz \right] w(u)\, du$$

$$= o(1).$$

This yields $T_3 = o(1)$. Combining these results, we obtain that

$$I_{skj} = \mu^{-1}(f_0)[w(k/s)/f_0(k/s)](1 + o(1)) = I_{sk}(1 + o(1)).$$

Now we can straightforwardly follow the proof of Theorem 1 in Efromovich (1989). This yields, for $k \in \{0, 1, \ldots, s-1\}$, that

$$\inf R_k \geq (s^{-2m} Q_{sk})^{1/(2m+1)} (nI_{sk})^{-2m/(2m+1)} P(1 + o(1)),$$

where the infimum is over all possible nonparametric estimates of $f$ considered in the theorem, and $P = (2m/2\pi(m+1))^{2m/(2m+1)}(2m+1)^{1/(2m+1)}$ is the Pinsker constant. Thus,

$$\inf \sum_{k=0}^{s-1} R_k \geq PQ^{1/(2m+1)} \left[ s^{-1} \mu(f_0) \sum_{k=0}^{s-1} f_0(k/s)/w(k/s) \right]^{2m/(2m+1)}$$

$$\times n^{-2m/(2m+1)}(1 + o(1))$$

$$= PQ^{1/(2m+1)} \left[ \mu(f_0) \sum_{k=0}^{s-1} \int_{k/s}^{(k+1)/s} (f_0(x)/w(x))\, dx \right]^{2m/(2m+1)}$$

$$\times n^{-2m/(2m+1)}(1 + o(1))$$

$$= PQ^{1/(2m+1)} \left( n^{-1} \mu(f_0) \int_0^1 (f_0(x)/w(x))\, dx \right)^{2m/(2m+1)} (1 + o(1)).$$

Theorem 1 is proved. $\square$

PROOF OF THEOREM 2. This proof follows along the lines of the proof of Theorem 1. Necessary changes are as follows. First, a new class $\mathcal{H}_s$ is introduced,

$$\mathcal{H}_s = \Big\{ f : f(x) = f_{0J_n}(x)$$

(4.8)
$$+ \left[ \sum_{k=0}^{s-1} f_k(x) - \sum_{k=1}^{s-1} \int_0^1 f_k(u) \left( 1 + \sum_{i=1}^{J_n} \varphi_i(u)\varphi_i(x) \right) du \right],$$

$$x \in [0,1], f_k \in \mathcal{H}_{sk}, f(x) = f_0(x), x \notin [0,1], f \geq 0 \Big\},$$

with no change in $\mathcal{H}_{sk}$ except for using $q_n Q$ in place of $Q$. Note that according to Assumption A the sequence $q_n$ decreases at most logarithmically and thus $J_n = o(1)J(k)$. Thus this class is defined correctly.

Second, we verify that for large $n$ the inclusion $\mathcal{H}_s \subset \mathcal{H}(m, Q, f_0, J_n)$ holds. Denote

$$\theta_j = \int_0^1 (f(u) - f_{0J_n}(u))\varphi_j(u)\, du, \qquad f \in \mathcal{H}_s.$$

Note that $\theta_j = 0$, $1 \le j \le J_n$, and thus the inclusion follows from the inequality

(4.9) $$\sum_{j > J_n} (\pi j)^{2m} \theta_j^2 \le q_n Q, \qquad f \in \mathcal{H}_s.$$

To prove (4.9), denote $\psi(u) = f(u) - f_{0J_n}(u)$. Because $\psi^{(s)}(0) = \psi^{(s)}(1) = 0$ for all odd $s < m$, using integration by parts implies [see Efromovich (1999), Section 2.2] $\theta_j^2 = (\pi j)^{-2m}[\int_0^1 \psi^{(m)}(u)\tilde{\varphi}_j(u)\, du]^2$, where $\tilde{\varphi}_j(u) = \varphi_j(u)$ for $m$ even and $\tilde{\varphi}_j(u) = \sqrt{2}\sin(\pi j x)$ for $m$ odd. This together with the Parseval identity implies, for both odd and even $m$, that

$$\sum_{j > J_n} (\pi j)^{2m} \theta_j^2 = \int_0^1 [\psi^{(m)}(u)]^2\, du.$$

Then following along the lines of (4.1) and (4.2) we get that

$$\int_0^1 [\psi^{(m)}(u)]^2 du \le q_n Q(1 - s^{-2}).$$

Inequality (4.9) is verified.

Finally note that in the estimation of $I_{skj}$, we get a new factor

$$\left[\varphi_{skj}(u)\phi_{sk}(u) - \int_{k/s}^{(k+1)/s} \varphi_{skj}(z)\phi_{sk}(z)\left(1 + \sum_{i=1}^{J_n} \varphi_i(z)\varphi_i(u)\right) dz\right]^2$$

in place of

$$\left[\varphi_{skj}(u)\phi_{sk}(u) - \int_{k/s}^{(k+1)/s} \varphi_{skj}(z)\phi_{sk}(z)\, dz\right]^2.$$

To evaluate this new factor, we write, similarly to (4.7),

$$\int_{k/s}^{(k+1)/s} \varphi_{skj}(z)\phi_{sk}(z)\left(1 + \sum_{i=1}^{J_n} \varphi_i(z)\varphi_i(u)\right) dz$$

$$= o(1)/\ln(n) + \int_{k/s}^{(k+1)/s} \varphi_{skj}(z)\left(1 + \sum_{i=1}^{J_n} \varphi_i(z)\varphi_i(u)\right) dz$$

$$= o(1)/\ln(n) + \sum_{i=1}^{J_n} \varphi_i(u) \int_{k/s}^{(k+1)/s} \varphi_{skj}(z)\varphi_i(z)\, dz.$$

Relationship (2.2.7) in Efromovich (1999) implies that

$$\left| \int_{k/s}^{(k+1)/s} \varphi_{skj}(z)\varphi_i(z)\,dz \right| \leq C j^{-2} i^2.$$

This inequality allows us to conclude that $I_{skj} = I_{sk}(1 + o(1))$, and then we finish the proof following along the lines of the proof of Theorem 1. $\square$

PROOF OF THEOREM 3. Here we are verifying the more complicated assertion (2.15). The assertion (2.14) is verified similarly and its proof is skipped.

The method of establishing sharp optimality of the Efromovich–Pinsker estimator is well developed and it consists of several steps. First of all, sharp optimality of a pseudoestimate is established. This is the step that should be verified for each particular problem. Then this estimate is mimicked by a so-called linear oracle that always performs better than the estimate. This step is easily verified. The third step is to show that the Efromovich–Pinsker blockwise oracle sharply mimics the linear oracle. For the particular blocks and thresholds considered in Section 2, this step is verified in Efromovich (1985) and it is well known. Finally, it should be shown that the Efromovich–Pinsker estimator mimics the Efromovich–Pinsker oracle. The validity of this step follows from Efromovich (1985, 2000).

Thus in what follows we verify steps 1 and 2. Consider a pseudoestimate

$$(4.10) \qquad \tilde{f}_n(u) = f_{0,J_n}(u) + \sum_{j=J_n+1}^{J^*} [1 - (j/J^*)^m]\hat{\theta}_j \varphi_j(u),$$

where $\hat{\theta}_j$ is defined in (2.12) and $J^*$ is the rounded up $[nd^{-1}(f,w)(2m+1) \times (m+1)q_n Q(2m(2\pi)^{2m})^{-1}]^{1/(2m+1)}$, $d(f,w) = \mu(f) \int_0^1 f(u)w^{-1}(u)\,du$. Then, for $f \in \mathcal{H}(m,Q,f_0,J_n)$,

$$E_f\left\{ \int_0^1 (\tilde{f}_n(u) - f(u))^2\,du \right\}$$

$$= \sum_{j=J_n+1}^{J^*} E_f\{[(1-(j/J^*)^m)\hat{\theta}_j - \theta_j]^2\} + \sum_{j>J^*} \theta_j^2$$

$$(4.11) \qquad = \sum_{j=J_n+1}^{J^*} [(1-(j/J^*)^m)^2 E_f\{(\hat{\theta}_j - \theta_j)^2\}$$

$$- 2(1-(j/J^*)^m)(j/J^*)^m \theta_j (E\{\hat{\theta}_j\} - \theta_j) + (j/J^*)^{2m}\theta_j^2]$$

$$+ \sum_{j>J^*} \theta_j^2.$$

Denote $\tilde{\theta}_j = n^{-1}\mu \sum_{l=1}^n \mathbb{1}(0 \leq Y_l \leq 1)w^{-1}(Y_l)\varphi_j(Y_l)$, that is, $\tilde{\theta}_j$ is the Cox empirical estimate $\hat{\theta}_j$ with $\hat{\mu}$ replaced by $\mu$. Then using the elementary identity

$$\hat{\mu} - \mu = \mu^2[(\mu^{-1} - \hat{\mu}^{-1}) - \hat{\mu}\mu^2(\mu^{-1} - \hat{\mu}^{-1})^2],$$

we get [in what follows $\mu = \mu(f)$]

$$\begin{aligned}
E_f\{\hat{\theta}_j\} &= E_f\{\tilde{\theta}_j + (\hat{\mu} - \mu)\mu^{-1}\tilde{\theta}_j\} \\
(4.12) \quad &= \theta_j + \mu^{-1}E_f\{(\hat{\mu} - \mu)\theta_j + (\hat{\mu} - \mu)(\tilde{\theta}_j - \theta_j)\} \\
&= \theta_j - \mu^3\theta_j E_f\{\hat{\mu}(\mu^{-1} - \hat{\mu}^{-1})^2\} + \mu^{-1}E_f\{(\hat{\mu} - \mu)(\tilde{\theta}_j - \theta_j)\}.
\end{aligned}$$

Also

$$\begin{aligned}
E_f\{(\hat{\theta}_j - \theta_j)^2\} &= E_f\{[(\tilde{\theta}_j - \theta_j) + (\hat{\mu} - \mu)\mu^{-1}\tilde{\theta}_j]^2\} \\
(4.13) \quad &= E_f\{(\tilde{\theta}_j - \theta_j)^2\} + 2\mu^{-1}E_f\{(\tilde{\theta}_j - \theta_j)(\hat{\mu} - \mu)\tilde{\theta}_j\} \\
&\quad + \mu^{-2}E\{(\hat{\mu} - \mu)^2\tilde{\theta}_j^2\}.
\end{aligned}$$

Using trigonometric relationships (3.1.7) and (3.1.8) in Efromovich (1999), we get

$$\begin{aligned}
(4.14) \quad & E_f\{(\tilde{\theta}_j - \theta_j)^2\} \\
&\leq n^{-1}\bigg[d(f,w) + \mu \int_0^1 f(u)w^{-1}(u)2^{-1/2}\varphi_{2j}(u)\,du\bigg],
\end{aligned}$$

$$\begin{aligned}
(4.15) \quad & E_f\{(\tilde{\theta}_j - \theta_j)(\hat{\mu} - \mu)\tilde{\theta}_j\} \\
&= \theta_j E_f\{(\tilde{\theta}_j - \theta_j)(\hat{\mu} - \mu)\} + E_f\{(\tilde{\theta}_j - \theta_j)^2(\hat{\mu} - \mu)\}
\end{aligned}$$

and

$$(4.16) \quad E_f\{(\hat{\mu} - \mu)^2\tilde{\theta}_j^2\} = \theta_j^2 E_f\{(\hat{\mu} - \mu)^2\} + E_f\{(\hat{\mu} - \mu)^2(\tilde{\theta}_j^2 - \theta_j^2)\}.$$

Then using the Cauchy–Schwarz inequality,

$$E_f\{(\hat{\mu} - \mu)^4\} < Cn^{-2}, \qquad E_f\{(\tilde{\theta}_j - \theta_j)^4\} < Cn^{-2},$$

we get

$$(4.17) \quad E_f\{(\hat{\theta}_j - \theta_j)^2\} \leq n^{-1}d(f,w) + n^{-1}\kappa_j + cn^{-3/2},$$

where $\sum_{j=1}^\infty \kappa_j^2 < C$.

Combining all these results in (4.11) we get

$$E_f\bigg\{\int_0^1 (\tilde{f}_n(u) - f(u))^2\bigg\}$$

$$
\text{(4.18)} \quad \leq \sum_{j=J_n+1}^{J^*} (1-(j/J^*)^m)^2 n^{-1} d(f,w)(1+o(1))
$$

$$
+ C n^{-1} \sum_{j=J_n+1}^{J^*} |\theta_j| + \sum_{j>J_n} (j/J^*)^{2m} \theta_j^2.
$$

Note that $\sum_{j=J_n}^{J^*} |\theta_j| = o(1)$ and $(\pi J^*)^{-2m} \sum_{j>J_n} (\pi j)^{2m} \theta_j^2 \leq (\pi J^*)^{-2m} q_n Q$ whenever $f \in \mathcal{H}(m,Q,f_0,J_n)$. Finally, plugging in $J^*$ and elementary calculations imply

$$
\text{(4.19)} \quad E_f \left\{ \int_0^1 (\tilde{f}_n(u) - f(u))^2 \right\} \leq [I_{fw} q_n^{-1/2m} n]^{-2m/(2m+1)} (1+o(1)).
$$

The first step in the proof is done. Then similarly to Section 7.4.5 in Efromovich (1999), we establish that the linear oracle

$$
\text{(4.20)} \quad f_{1n}(u) = f_{0J_n}(u) + \sum_{j=J_n+1}^{n^{1/3}} \theta_j^2 (\theta_j^2 + \hat{d} n^{-1})^{-1} \hat{\theta}_j \varphi_j(u)
$$

dominates the pseudoestimate (4.10).

Finally note that the elementary relationship

$$
E_f \left\{ \sum_{j=0}^{J_n} (\hat{\theta}_j - \theta_j)^2 \right\} \leq C \ln(n) n^{-1}
$$

holds. This allows us to follow along the lines of Efromovich (1985, 2000) and to verify the above-described last two steps in the proof. □

**Acknowledgment.** Comments of a referee are appreciated.

Department of Mathematics and Statistics
The University of New Mexico
Albuquerque, New Mexico 87131
USA
e-mail: efrom@math.unm.edu